\documentclass[11pt,reqno]{amsart}
\usepackage{amssymb,amscd,amsbsy,mathrsfs}
\usepackage{mathdots}
\newcommand{\lb}{\linebreak}

\renewcommand{\d}{\delta}
\newcommand{\e}{\varepsilon}

\newcommand{\s}{\sigma}

\newcommand{\f}{\varphi}

\newcommand{\C}{{\Bbb C}}

\newcommand{\R}{{\Bbb R}}
\newcommand{\Z}{{\Bbb Z}}

\newcommand{\bS}{{\boldsymbol S}}

\newcommand{\rf}[1]{(\ref{#1})}

\newcommand{\df}{\stackrel{\mathrm{def}}{=}}

\newcommand{\const}{\operatorname{const}}

\newcommand{\eeq}{\end{equation}}
\newcommand{\beq}{\begin{equation}}
\newcommand{\bay}{\begin{eqnarray}}
\newcommand{\ba}{\begin{align*}}
\newcommand{\ea}{\end{align*}}
\newcommand{\ey}{\end{eqnarray}}
\newcommand{\bey}{\begin{eqnarray*}}
\newcommand{\eey}{\end{eqnarray*}}

\newcommand{\be}{\infty}

\newcommand{\bl}{\blacksquare}

\newcommand{\Pf}{{\bf Proof. }}

\newtheorem{thm}{\hspace{\parindent}Theorem}[section]

\newtheorem{cor}[thm]{\hspace{\parindent}Corollary}
\newtheorem{lem}[thm]{\hspace{\parindent}Lemma}

\theoremstyle{remark}

\newtheorem*{rem*}{Remark}

\newcommand\OL{{\rm OL}}

\newcommand{\fF}{{\frak F}}

\newcommand\Li{{\rm Lip}}

\newcommand{\Bbbone}{{\rm{1\mathchoice{\kern-0.25em}{\kern-0.25em}{\kern-0.2em}{\kern-0.2em}I}}}

\begin{document}

\numberwithin{equation}{section}

\numberwithin{equation}{section}

\title{Functions of compact operators under trace class perturbations}
\author{A.B. Aleksandrov and V.V. Peller}
\thanks{The research  is supported by 
Russian Science Foundation [grant number 23-11-00153]. 
}


\

\begin{abstract}
The paper studies the problem, for which continuous functions $f$ on the real line $\R$,
the difference of the functions $f(B)-f(A)$ of self-adjoint operators $A$ and $B$ with trace class difference
must also be of trace class. The main result of the paper shows that this happens if and only if the function $f$ is operator Lipschitz on a neighbourhood of zero.
\end{abstract} 

\maketitle

\setcounter{section}{0}
\section{\bf Introduction}
\setcounter{equation}{0}
\label{In}

\

In this paper we continue studying the behaviour of functions of self-adjoint operators under their perturbation and we refer the readerto earlier papers \cite{Pe1} and \cite{Pe2} as well as to the survey
\cite{APOL}. 

We sre interested in the problem of when trace class perturbations of an operator lead to trace class increments of the functions of the operator. This property is closely related to the property of operator Lipschitzness. 

In this paper we consider this problem for functions of compact self-adjoint operators. To be more precise, 
we will be interested in the problem of describing the class of functions $f$ on the real line $\R$ that possess that property:

\medskip

{\it Let $A$ and $B$ be compact self-adjoint operators such that $A-B\in\bS_1$.
Then $f(A)-f(B)\in\bS_1$}.

\medskip

We denote by $\bS_1$ the class of trace class operators on Hilbert space.
We refer the reader to the books \cite{GK} and \cite{BS} for detailed information on trace class operators and other Schatten-von Neumann classes.

The main result of the paper is given in \S\:\ref{osnova}. We show that the class of functions in question consists of functions $f$ on $\R$ such that the restriction of $f$ to a certain interval $[-\e,\e]$, $\e>0$, is an operator Lipschitz function on that interval.

In \S\:\ref{opr} we give necessary definitions while in the concluding section \S\:\ref{normal'no} we briefly discuss a similar problem for functions of compact normal operators with trace class difference.

\

\section{\bf Basic definitions}
\setcounter{equation}{0}
\label{opr}

\

As we have already mentioned in the introduction, we are going to investigate in which cases a function defined on the real line or on a subset of it has the property that trace class perturbations lead to trace class increments of functions of such operators. 

\medskip

{\bf Definition 1.} Let $\fF$ be a closed subset of the real line and let $f$ be a continuous complex-valued function on $\fF$. We say that $f$ {\it respects trace class}, if the facts that
$A$ and $B$ are self-adjoint (not necessarily bounded!) operators on Hilbert space with spectra in 
$\fF$ and $B-A\in\bS_1$ imply that $f(B)-f(A)\in\bS_1$.

\medskip

It turns out that this property is closely related to the notion of operator Lipschitz functions and the notion of trace class Lipschitz functions.

Note here that in Farfarovskaya's paper \cite{F2} it was shown that not all Lipschitz functions on the real line respect trace class. Note also here that in the papers \cite{Pe1} and \cite{Pe2} it was shown that the functions of the Besov class $B_{\be,1}^1$ on the real line respect trace.

It turns out that this property that this property is equivalent to two other properties that we are going to discuss. 

\medskip

{\bf Definition 2.}  Let $\fF$ be a closed subset of the real line. We say that a complex-valued function $f$ on 
$\fF$ is {\it is an operator Lipschitz function} on $\fF$ if the inequality
\bay
\label{OLf}
\|f(A)-f(B)\|\le\const\|A-B\|
\ey
holds for arbitrary self-adjoint operators $A$ and $B$ on Hilbert space whose spactra $\s(A)$ and $\s(B)$ are subset of $\fF$.  Denote by $\OL(\fF)$ the space of all operator Lipschitz functions on $\fF$. The norm (to be more precise, the seminorm) of an operator Lipschitz function $f$ is, by definition, the number
$$
\|f\|_{\OL(\fF)}\df\sup\frac{\|f(A)-f(B)\|}{\|A-B\|},
$$
where the supremum is taken over all self-adjoint operators $A$ and $B$ such that
$A-B$ is a bounded operator and $A\ne B$.

\medskip

Note also that in this definition of $\OL(\fF)$ we can consider only bounded operators as well as not necessarily bounded operators; in both cases we get the same class of functions. Also in the definition of the norm
$\|\cdot\|_{\OL(\fF)}$ it does not matter whether we consider only bounded operators or we allow unbounded operators. 

The class of operator Lipschitz functions on $\R$ has rather amusing properties. For example, 
the functions in $\OL(\R)$ must be differentiable at every point of the real line which was proved in \cite{W}. On the other hand, operator Lipschitz functions on $\R$ do not have to be continuously differentiable  which was observed in \cite{KS}.

We refer the reader to the survey \cite{APOL} in which such issues are discussed in detail.

Recall also for completeness the definition of LIpschitz functions.
Let $\fF$ be a closed subset of the real line.
Denote $\Li(\fF)$ by the space of (complex-valued) Lipschitz functions on $\fF$ that is
the set of functions $f:\fF\to\C$ such that
\bay
\label{Lif}
|f(x)-f(y)|\le\const|x-y|
\ey
for all $x,y\in\fF$. The Lipschitz norm (to be more precise, the seminorm) $\|f\|_{\Li(\fF)}$ is defined as the smallest constant satisfying inequality
\rf{Lif} for all $x,y\in\fF$.

It is easy to see that each operator Lipschitz function on $\fF$ is Lipschitz and
$$
\|f\|_{\Li(\fF)}\le\|f\|_{\OL(\fF)},\quad f\in \OL(\fF).
$$

It turned out, however, that the converse is false which was shown for the first time by Farforovskaya's  in  \cite{F}.

Let us proceed now to the property of trace class Lipschitzness.

\medskip

{\bf Definition 3.} Let $\fF$ be a closed subset of the real line. We say that a complex-valued continuous function  $f$ on $\fF$ is {\it trace class Lipschitz} on $\fF$ if 
\bay
\label{OLS}
\|f(A)-f(B)\|_{\bS_1}\le\const\|A-B\|_{\bS_1}
\ey
for arbitrary self-adjoint operators $A$ and $B$ with spectra in $\fF$ such that  
$A-B\in\bS_1$.

\medskip

Herewith the trace class Lipschitz norm $\|f\|_{\OL_{\bS_1}(\fF)}$ is defined as the smallest constant, for which 
inequality \rf{OLS} holds for all self-adjoint operators $A$ and $B$ such that $\s(A), \s(B)\subset\fF$
and $A-B\in\bS_1$.

Theorem 3.6.3 of the paper \cite{APOL} says that for a function $f$ of class $\OL(\fF)$ the following inequalities hold:
\bay
\label{OS}
\frac12\|f\|_{\OL(\fF)}\le\|f\|_{\OL_{\bS_1}(\fF)}\le 2\|f\|_{\OL(\fF)}.
\ey
In the same survey \cite{APOL} in Theorem 3.6.4 it is proved that $\|f\|_{\OL_{\bS_1}(\fF)}=\|f\|_{\OL(\fF)}$ if
$\fF$ is a perfect subset of the real line.

It can be reduced fairly easily from \rf{OS} that for a continuous function $f$ on a closed subset $\fF$ of the real line, the following statements are equivalent:

(a) $f$ is an operator Lipschitz function on $\fF$;

(b) the function $f$ respects trace class;

(c) $f$ is a trace class Lipschitz function on $\fF$.

\medskip

This assertion was stated explicitly in the survey \cite{APOL} in the case when $\fF=\R$.

\

\section{\bf The main result}
\setcounter{equation}{0}
\label{osnova}

\

In this section we obtain the main result of the paper, which gives a description of the class of continuous functions $f$ on the real line, for which the assumptions that $A$ and $B$ are compact self-adjoint operators and $B-A\in\bS_1$ imply that $f(B)-f(A)\in\bS_1$.

To prove the result, we use the observation made in the survey 
\cite{APOL}, which says that to establish the operator Lipschitzness as well as the trace class Lipschitzness, 
it suffices to verify those properties on operators on finite-dimensional spaces.

\begin{lem}
\label{yadl}
Let $f$ be a continuous function on a closed set $\fF$, $\fF\subset\R$.
Suppose that $\|f\|_{\OL(\fF)}>m>0$. Then there exist self-adjoint operators 
$A$ and $B$ on a finite-dimensional Hilbert space such that
$\s(A),\s(B)\subset \fF$ and
$$
\|f(B)-f(A)\|_{\bS_1}>\frac m2\|B-A\|_{\bS_1}>0.
$$
\end{lem}

\Pf It is easy to verify that there exists a finite subset $\fF_0$ of  $\fF$ such that
$\|f\|_{\OL(\fF_0)}>m>0$, see, for example, dormula (3.1.8) in the survey \cite{APOL}.
Then there exist self-adjoint operators $A$ and $B$ on a finite-dimensional Hilbert space such that
$$
\|f(B)-f(A)\|_{\bS_1}=\|f\|_{\OL_{\bS_1}(\fF_0)}\|B-A\|_{\bS_1}>0.
$$
It remains to observe that $2\|f\|_{\OL_{\bS_1(\fF_0)}}\ge\|f\|_{\OL(\fF_0)}>m$. $\bl$

\begin{cor}
\label{yadc}
Let $f$ be a continuous function defined on a nondegenerate interval  $[a,b]$, $-\be<a<b<\be$.
Suppose that $\|f\|_{\OL([a,b])}>m$. Then there exist self-adjoint operators $A$ and $B$ on a finite-dimensional Hilbert space such that
$\s(A),\s(B)\subset [a,b]$,
$$
\frac{\|f(B)-f(A)\|_{\bS_1}}{\|B-A\|_{\bS_1}}>\frac m2 \quad \text {and} \quad \|f(B)-f(A)\|_{\bS_1}<1.
$$
\end{cor}

\Pf Clearly, there exist self-adjoint operators $A$ and $B$ on a finite-dimensional Hilbert space such that
$\s(A),\s(B)\subset [a,b]$ and
$$
\frac{\|f(B)-f(A)\|_{\bS_1}}{\|B-A\|_{\bS_1}}>\frac m2.
$$
Suppose that $\|f(B)-f(A)\|_{\bS_1}\ge1$.
Consider the function $\f(t)\df f((1-t)A+tB)$ ($t\in[0,1]$). By the Cantor theorem, this function is uniformly continuous.
Consequently, there exists a positive integer $n$ such that 
$\big\|\f\big(\frac{k+1}n\big)-\f\big(\frac kn\big)\big\|_{\bS_1}<1$
for all $k$ in $[0,n)\cap\Z$. Clearly, 
$$
\Big\|\f\Big(\frac{k+1}n\Big)-\f\Big(\frac kn\Big)\Big\|_{\bS_1}\ge\frac1n\|\f(1)-\f(0)\|_{\bS_1}\ge\frac1n\|f(B)-f(A)\|_{\bS_1}
$$
for a certain $k$ in $[0,n)\cap\Z$. Then
$$
\frac{\|f(A+\frac{k+1}n(B-A))-f(A+\frac{k}n(B-A))\|_{\bS_1}}{\frac1n\|B-A\|_{\bS_1}}
\ge\frac{\|f(B)-f(A)\|_{\bS_1}}{\|B-A\|_{\bS_1}}>\frac m2. \quad \bl
$$


\begin{cor}
\label{yadc2}
Let $f$ be a contibuous function on a closed nondegenerate interval $P$.
Suppose that $\|f\|_{\OL(P)}>m>0$. Then there exist self-adjoint operators
$A$ and $B$ on a finite-dimensional Hilbert space such that
$\s(A),\s(B)\subset P$,
$$
\frac{\|f(B)-f(A)\|_{\bS_1}}{\|B-A\|_{\bS_1}}>\frac m2 \quad \text {and} \quad \frac12\le\|f(B)-f(A)\|_{\bS_1}\le1.
$$
\end{cor}

\Pf We may assume that the interval $P$ is bounded. Let $A$ and $B$ be operators in the conclusion of Corollary \ref{yadc}.
Put $N=\big[\|f(B)-f(A)\|^{-1}_{\bS_1}\big]$. Here for a real number $x$, we use the notation $[x]$ for the largest integer that does not exceed $x$. 

We define the operators $A_N$ and $B_N$ by
$$
A_N=\bigoplus_{j=1}^{N}A\quad\mbox{and}\quad B_N=\bigoplus_{j=1}^{N}B.
$$
Then
$$
\frac12\le N\|f(B)-f(A)\|_{\bS_1}=\|f(B_N)-f(A_N)\|_{\bS_1}\le 1.
$$
It remains to observe that
$$
\frac{\|f(B_N)-f(A_N)\|_{\bS_1}}{\|B_N-A_N\|_{\bS_1}}=\frac{\|f(B)-f(A)\|_{\bS_1}}{\|B-A\|_{\bS_1}}>\frac m2
$$ 
and
$$
\frac12\le\|f(B_N)-f(A_N)\|_{\bS_1}\le1. \quad \bl
$$

Before we proceed to the main result of the paper, we consider a simpler problem for functions of {\it commuting} 
compact self-adjoint operators with trace class difference.

\begin{lem}
\label{yadl2}
Let $f$ be a function on the real line. Then the following statements are equivalent:

{\em(a)} the implication 
$$
A-B\in\bS_1\Longrightarrow f(A)-f(B)\in\bS_1
$$
is valid for arbitrary commuting compact self-adjoint operators $A$ and $B$;

{\em(b)} there exists a positive number $\d>0$ such that $f\big|[-\d,\d]\in\Li([-\d,\d])$.
\end{lem}

\Pf Let us observe that since compact self-adjoint operators $A$ and $B$  commute, it follows that there exists an orthonormal basis of vectors that are eigenvectors of both operators $A$, and $B$.
Thus, condition (a) means that for sequences  $\{t_k\}_{k=1}^\be$ and $\{s_k\}_{k=1}^\be$ that tend to zero, 
the following implication is valid:
\bay
\label{CK}
\sum_{k=1}^\be|t_k-s_k|<\be \quad\Longrightarrow\quad \sum_{k=1}^\be|f(t_k)-f(s_k)|<\be.
\ey
Therefore, the implication (b)$\Rightarrow$(a) is trivial. 
To show that (a) implies (b), we observe that the implication \rf{CK} can be rewritten in the following way 
\bay
\label{CK1}
\sum_{k=1}^\be n_k|t_k-s_k|<\be \quad\Longrightarrow\quad \sum_{k=1}^\be n_k|f(t_k)-f(s_k)|<\be,
\ey
where $\{n_k\}_{k=1}$ is an arbitrary sequence of positive integers. 
Suppose that (b) is not valid. Then there exist sequences
$\{t_k\}_{k=1}^\be$ and $\{s_k\}_{k=1}^\be$ that tend to zero and such that  \lb$0<|t_k-s_k|<2^{-k}$ and $\frac{|f(t_k)-f(s_k)|}{|t_k-s_k|}>2^k$
for all $k$ in $\Bbb N$. It remains to find a sequence $\{n_k\}_{k=1}^\be$ of positive integers
such that
$$
\sum_{k=1}^\be n_k|t_k-s_k|<\be
$$
and
$$
\sum_{k=1}^\be n_k|f(t_k)-f(s_k)|=\be.
$$

Put $n_k=\big[|f(t_k)-f(s_k)|^{-1}\big]+1$. Then 
$$
\sum_{k=1}^\be n_k|f(t_k)-f(s_k)|=\be
$$
because
$n_k|f(t_k)-f(s_k)|\ge1$ for all $k$ in $\Bbb N$. It remains to prove that 
$$
\sum_{k=1}^\be n_k|t_k-s_k|<\be.
$$
To prove this, it remains to observe that
$$
n_k|t_k-s_k|\le|t_k-s_k|+\frac{|t_k-s_k|}{|f(t_k)-f(s_k)|}<2^{1-k}.\quad\bl
$$

\medskip

Let us proceed now to the main result of the paper.

\begin{thm}
\label{yad}
Let $f$ be a function on the real line.
Then the following statements are equivalent:

{\em(a)} the implication 
$$
A-B\in\bS_1\Longrightarrow f(A)-f(B)\in\bS_1
$$
holds for arbitrary compact self-adjoint operators $A$ and $B$;

{\em(b)} there exists a positive number $\d$ such that $f\big|[-\d,\d]\in\OL([-\d,\d])$.
\end{thm}

\Pf Let us first prove that (b)$\Rightarrow$(a).  We have to show that
$$
A-B\in\bS_1\Longrightarrow f(A)-f(B)\in\bS_1.
$$
Obviously, we may assume that $f(0)=0$.
Consider the operators $A_\d\df A\chi_{[-\d,\d]}(A)$ and
$B_\d\df B\chi_{[-\d,\d]}(B)$.  First of all, we observe that $A-A_\d$, $B-B_\d$, $f(A)-f(A_\d)$ and $f(B)-f(B_\d)$ are finite rank operators, and so they are of trace class. Thus, it suffices to show that 
$$
A_\d-B_\d\in\bS_1\Longrightarrow f(A_\d)-f(B_\d)\in\bS_1.
$$
However, this is clear because $\s(A_\d),\s(B_\d)\subset(-\d,\d)$ and $f\in\OL(-\d,\d)$.

Let us now show that (a)$\Rightarrow$(b). Assume the contrary, that is the restriction of $f$ to any interval 
$[-\d,\d]$ is ot operator Lipschitz on $[-\d,\d]$ for any positive number $\d$.

It follows from Lemma \ref{yadl2} that $f\in\Li([-\d,\d])$ for a certain positive number $\d$.
Consider a sequence $\{\d_n\}_{n=1}^\infty$ of points of $(0,\d)$ that tends to zero.
It follows from Corollary \ref{yadc2} that for each positive integer $n$ there exist self-adjoint operators 
$A_n$ and $B_n$ on a finite-dimensional HIlbert space such that 
$\|A_n\|,\|B_n\|<\d_n$,
$$
\frac{\|f(B_n)-f(A_n)\|_{\bS_1}}{\|B_n-A_n\|_{\bS_1}}>2^n \quad \text {and} \quad \frac12\le\|f(B_n)-f(A_n)\|_{\bS_1}\le1.
$$
Put $A=\bigoplus_{j=1}^{\be}A_n$ and $B=\bigoplus_{j=1}^{\be}B_n$. Clearly, $A$ and $B$ are compact self-adjoint operators such that 
$$
\|B-A\|_{\bS_1}=\sum_{n=1}^\be\|B_n-A_n\|_{\bS_1}<\be
$$
and
$$
\|f(B)-f(A)\|_{\bS_1}=\sum_{n=1}^\be\|f(B_n)-f(A_n)\|_{\bS_1}=\be. \quad \bl
$$

\begin{cor}
Suppose that $f$ is a continuous function on the real line that satisfies condition {\em(a)} in the statement of Theorem {\em\ref{yad}}. Then the function $f$ must be differentiable in a neighbourhood of the origin. 
\end{cor}

\begin{cor}
There exist compact self-adjoint operators $A$ and $B$ on Hilbert space such that 
$$
A-B\in\bS_1\quad\mbox{but}\quad|A|-|B|\not\in\bS_1.
$$
\end{cor}

\

\section{\bf Normal operators with trace class difference}
\setcounter{equation}{0}
\label{normal'no}

\

A natural generalization of the problem solved in \S\:\ref{osnova} would be a similar problem for functions of compact operators with trace class difference. 

Let $M$  and $N$ be compact normal operators such that $N-M\in\bS_1$. The problem is to describe the class of continuous functions $f$ on he complex plane $\C$ such that under such conditions the difference $f(N)-f(M)$ 
of functions of the operator $N$ and $M$ must also be of trace class.

Unfortunately, our method of the proof of the main result does not extend to the case of functions of normal operators. In the case of continuous functions defined on closed subsets of the complex plane it is also possible to define the classes of operator Lipschitz functions and trace class Lipschitz functions. However, unlike in the case of functions on subsets of the real line the authors still do not know whether these two classes have to coincide. Yet, there are also other problems. In particular, it is not quite clear whether the condition of operator Lipschitzness can be verified only for operators on finite-dimensional Hilbert spaces. 

%
%
%
%
%
%

\

\
 
%
%

\begin{footnotesize}
 
\noindent
\begin{tabular}{p{8cm}p{15cm}}
A.B. Aleksandrov & V.V. Peller \\
St.Petersburg State University & St.Petersburg State University \\
Universitetskaya nab., 7/9  & Universitetskaya nab., 7/9\\
199034 St.Petersburg, Russia & 199034 St.Petersburg, Russia \\
\\

St.Petersburg Department &St.Petersburg Department\\
Steklov Institute of Mathematics  &Steklov Institute of Mathematics  \\
Russian Academy of Sciences  & Russian Academy of Sciences \\
Fontanka 27, 191023 St.Petersburg &Fontanka 27, 191023 St.Petersburg\\
Russia&Russia\\
email: alex@pdmi.ras.ru& email: peller@math.msu.edu
\end{tabular}
\end{footnotesize}

\end{document}